\documentclass{amsbook}

\usepackage[T1]{fontenc}
\usepackage{ae,aecompl}

\usepackage{amssymb,amscd}
\usepackage{url}
\usepackage{amsthm,amsfonts,amsmath,amsxtra}
\usepackage{graphicx}

\usepackage{tikz}

\numberwithin{equation}{section}
\numberwithin{table}{section}
\numberwithin{figure}{section}
\numberwithin{section}{chapter}

\newtheoremstyle{bold}
{.5\baselineskip}{.5\baselineskip}{\itshape}{}{\bfseries}{.}{.5em}{}
\newtheoremstyle{shy}
{.5\baselineskip}{.5\baselineskip}{}{}{\bfseries}{.}{.5em}{}

\makeatletter
\def\@captionfont{\small}

\def\mychapter{%
  \if@openright\cleardoublepage\else\clearpage\fi
 \thispagestyle{empty}\global\@topnum\z@
  \@afterindenttrue \secdef\@mychapter\@schapter}

\def\@mychapter[#1]#2#3{\refstepcounter{chapter}%
  \ifnum\c@secnumdepth<\z@ \let\@secnumber\@empty
  \else \let\@secnumber\thechapter \fi
  \typeout{\chaptername\space\@secnumber}%
  \def\@toclevel{0}%
  \ifx\chaptername\appendixname \@tocwriteb\tocappendix{chapter}{#2\\ \scshape #3}%
  \else \@tocwriteb\tocchapter{chapter}{#2\\ \scshape #3}\fi
  \chaptermark{#1}%
  \addtocontents{lof}{\protect\addvspace{10\p@}}%
  \addtocontents{lot}{\protect\addvspace{10\p@}}%
  \@mymakechapterhead{#2}{#3}\@afterheading}

\def\@mymakechapterhead#1#2{\global\topskip 7.5pc\relax
  \begingroup
  \fontsize{\@xivpt}{18}\bfseries\centering
    \ifnum\c@secnumdepth>\m@ne
      \leavevmode \hskip-\leftskip
      \rlap{\vbox to\z@{\vss
          \centerline{\normalsize\mdseries
              \uppercase\@xp{\chaptername\ \thechapter}}
          \vskip 3pc}}\hskip\leftskip\fi
     #1\par \vskip 1pc
     \Large\mdseries\scshape\centering
     #2\par \endgroup
  \skip@34\p@ \advance\skip@-\normalbaselineskip
  \vskip\skip@ }

\def\section{\@startsection{section}{1}%
  \z@{.9\linespacing\@plus\linespacing}{.5\linespacing}%
  {\large\bfseries\boldmath\centering}}

\def\subsection{\@startsection{subsection}{2}%
  \z@{.7\linespacing\@plus\linespacing}{.5\linespacing}%
  {\normalfont\scshape\centering}}

\def\theindex{\@restonecoltrue\if@twocolumn\@restonecolfalse\fi
  \columnseprule\z@ \columnsep 35\p@
  \@indextitlestyle
  \thispagestyle{empty}%
  \let\item\@idxitem
  \parindent\z@  \parskip\z@\@plus.3\p@\relax
  \raggedright
  \hyphenpenalty\@M
  \footnotesize}

\renewcommand{\@bibtitlestyle}{%
  \@xp\section\@xp*\@xp{\bibname}%
}
\renewcommand{\tocchapter}[3]{%
  \indentlabel{\@ifnotempty{#2}{\ignorespaces#1 #2.\quad}}#3}

\renewcommand{\tocsection}[3]{%
  \indentlabel{\@ifnotempty{#2}{\makebox[3.2em][l]{\ignorespaces#1 #2.}}}#3}

\makeatother

\renewcommand{\tocappendix}[3]{%
  \indentlabel{#1.\quad}#3}
\renewcommand{\tocappendix}[3]{%
  \indentlabel{\makebox[5.7em][l]{\ignorespaces#1.}}#3}

\renewcommand{\bibname}{References}
\renewcommand{\ge}{\geqslant}
\renewcommand{\geq}{\geqslant}

\renewcommand{\leq}{\leqslant}

\theoremstyle{bold}
\newtheorem{theorem}{Theorem}[section]
\newtheorem{proposition}[theorem]{Proposition}

\newtheorem{corollary}[theorem]{Corollary}

\theoremstyle{shy}
\newtheorem{definition}[theorem]{Definition}
\newtheorem{remark}[theorem]{Remark}

\newcommand{\EE}{\mathbb{E}}

\newcommand{\NN}{\mathbb{N}}

\newcommand{\PP}{\mathbb{P}\ts}

\newcommand{\RR}{\mathbb{R}}

\newcommand{\dd}{\ts\mathrm{d}\ts}

\newcommand{\ts}{\hspace{0.5pt}}

\newcommand{\bz}{{\overline{z}}}
\newcommand{\br}{{\overline{r}}}

\newcommand{\bmu}{{\overline{\mu}}}
\newcommand{\bLambda}{{\overline{\Lambda}}}

\begin{document}


\begin{center}
{\huge\bf Population genetic models of dormancy}\\[1cm]

{\Large\sc Jochen~Blath and Noemi~Kurt}\\[1cm]
\end{center}

{\small 
In the present article, we investigate the effects of dormancy on an abstract population genetic level. We first provide a short review of seed bank models in population genetics, and the role of dormancy for the interplay of evolutionary forces in general, before we discuss two recent paradigmatic models, referring to spontaneous resp.\ simultaneous switching of individuals between the active and the dormant state. We show that both mechanisms give rise to non-trivial mathematical objects, namely the (continuous) {\em seed bank diffusion} and the {\em seed bank diffusion with jumps}, as well as their dual processes, the {\em seed bank coalescent} and the {\em seed bank coalescent with simultaneous switching}. }


\index{seed bank}
\index{coalescent}
\index{Kingman coalescent}
\index{dormancy}

\section[Introduction]{Introduction}

Recently, the phenomenon of individual \emph{dormancy} has attracted significant attention in a population genetic context (cf.\ e.g.\ \cite{JBNK-LJ11}, which inspired much of the presented research, and \cite{JBNK-SL18} for a recent systematic overview providing further references). The term `dormancy' here refers to the ability of an organism to enter a reversible state of low or zero metabolic activity, in which it does not reproduce, but simply persists, unaffected by other forces, for potentially extended periods of time.
Dormancy seems to be a wide-spread and important evolutionary trait that has been developed by many species in many different guises across the tree of life. It is often seen as a bet-hedging strategy that allows organisms to persist through unfavourable environmental conditions and leads to a seed bank, in which genotypic and phenotypic variability can be stored for extended periods of time. For example, many microbial populations maintain a reservoir of dormant individuals, and in fact it seems that large fractions of microbial populations are in a dormant state at any given time, \cite{JBNK-LJ11}. Seed banks have also been described and analysed with the help of population genetic quantities in other species (as an example, we mention populations of wild tomatoes investigated in \cite{JBNK-T11}). 
It is plausible to assume that dormancy affects classical evolutionary forces such as genetic drift, mutation and selection in substantial ways that should generally increase genetic variability. A systematic discussion of the effect of the presence of seed banks on other evolutionary forces can be found in  \cite{JBNK-SL18}, see also \cite{JBNK-L90, JBNK-V04}. It thus might seem justified to discuss whether the effect of dormancy itself should be considered as an evolutionary force in its own right. 

However, a systematic approach to the mathematical modelling of dormancy and seed banks in population genetics has only started relatively recently. 
A first sophisticated mathematical model in this direction has been introduced by Kaj, Krone and Lascoux \cite{JBNK-KKL01}. Here, the authors consider a variant of a classical Wright-Fisher model, say of size $N$, where offspring individuals do not necessarily ``choose'' their parent from the previous generation (as in the classical Wright-Fisher model), but from individuals alive an independent random number of generations ago (where the random number of generations is assumed to be bounded by some constant $m>0$). This situation is then interpreted as the offspring of the parental particle staying in a seed bank for the corresponding amount of generations. As usual in theoretical population genetics, the limit of large populations ($N \to \infty$) is studied on the `evolutionary timescale' measured in units of population size $N$ (see e.g.\ \cite{JBNK-W09} for an overview and further references).  If the bound $m$ on the number of generations in the above model is finite and in particular independent of $N$, then the time spent in the seed bank is short compared to the evolutionary timescale, which leads to what is sometimes called a \emph{weak} seed bank effect. Not unexpectedly, the genealogy that is obtained after the classical evolutionary rescaling is a constant time change of the Kingman coalescent, and therefore typical patterns of genetic diversity like the (normalized) site frequency spectrum remain qualitatively unchanged in this model. However, estimates of effective population size resp.\ coalescent mutation rate are affected, see e.g.\ \cite{JBNK-T11}. The mechanism has been extended to a potentially unbounded (but with finite mean) time in the seed bank \cite{JBNK-BGKS12} and to incorporate selection \cite{JBNK-MKTZ17}, see also \cite{JBNK-ZT12}. It appears that this model is tailored to populations, in which the time of latency in the seed bank is in the range of few generations (compared to the total population size), as e.g.\ in many plants.

\index{site frequency spectrum}
\index{latency}
\index{Wright-Fisher model}
\index{migration}
\index{island model}

\begin{figure}[t]
\label{fig:Scheme}

\begin{center}
\setlength{\unitlength}{3cm}

\begin{tikzpicture}

\thinlines

\thicklines

\draw[->] (-1.5,1) .. controls (0,1.4) .. (1.5,1);
\draw[<-] (-1.5,-1) .. controls (0,-1.4) .. (1.5,-1);
\draw[very thick] (-2.5,0) circle [radius=1];
\draw[very thick] (2.5,0) circle [radius=1];
\draw[<-] (-3.2,1) arc (0:270:4mm);

\put(-.97,-.02){$\rm {active \atop cells}$}
\put(0.65,-.02){$\rm {dormant\atop cells}$}
\put(-.5,.51){$\rm Resuscitation \, \, rate\, \, \it cK$}
\put(-.46,-.6){$\rm Dormancy \, \, rate\, \, \it cK$}
\put(-1.5,0.51){$\rm Reproduction$}
\end{tikzpicture}
\medskip
\medskip
\medskip

\caption{\label{JBNK-fig:schema} Schematic representation of a population with seed bank corresponding to two `islands' as in \cite{JBNK-LJ11}}
\end{center}

\end{figure}
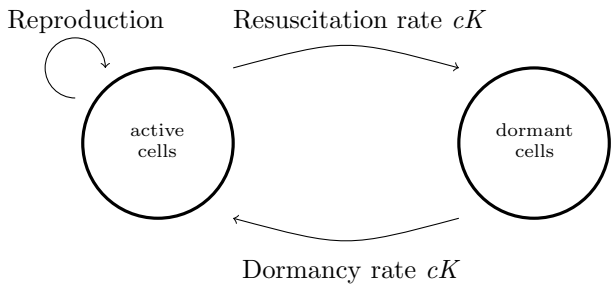

A different class of seed bank models has recently been constructed in the spirit of the Wright-Fisher model with two islands: Here, entering and leaving a dormant state is considered as `migration' between two `islands' (comprising of the active and the dormant population, cf. \cite{JBNK-BEGK15, JBNK-BGKW16}). This setup seems to fit to the case of bacterial communities and was suggested in \cite{JBNK-LJ11}, see Figure \ref{JBNK-fig:schema}. Here, a positive fraction of the population (of order $N$) will be in a dormant state, and the latency times spent in the dormant state will also necessarily be of order $N$, which is different from the assumption of the model \cite{JBNK-KKL01}. In such a situation, we speak of a {\em strong} seed bank effect. 


The presence of a seed bank of the above type will indeed drastically change the behaviour of the scaling limit of the population model as well as its genealogical process, and will produce rather unique patterns of genetic variability. This recent model, the corresponding ancestral process and its properties will be presented and discussed in Section \ref{JBNK-sect:def}.

Note that modeling seed banks in such a `migration' set-up implicitly assumes that individuals switch independently from active to dormant and vice versa. This corresponds to \emph{spontaneous switching} of bacteria as discussed in \cite{JBNK-LJ11}, and should be appropriate for populations in `stable' environments. However, often one observes that initiation of or resuscitation from dormancy is triggered by environmental cues, and in such situations one will see many individuals switch state at the same time, so that the independence assumption of classical migration is violated. This behaviour corresponds to the notion of \emph{responsive switching} from \cite{JBNK-LJ11}. We thus also present and discuss a very recent model which includes such simultaneous switching events between states in Section \ref{JBNK-sect:sim} based on the preprint \cite{JBNK-BGKW19}, introducing the \emph{seed bank coalescent with simultaneous switching}. Of course, it is also possible to include both switching mechanisms in a joint model (and in fact we will consider them jointly below in Section \ref{JBNK-sect:sim}). The resulting genealogies display a rich behaviour, depending on the choice of parameters governing the simultaneous versus spontaneous switching mechanisms. Note that the above models also allow the derivation of inference methods, with some results on the spontaneous switching case to be found in \cite{JBNK-BEGK15} and the preprint \cite{JBNK-BBKWB18+}. For the simultaneous switching case, this is work in progress. 

An interesting variant of the above switching regimes is presented in \cite{JBNK-BC+16}, \cite{JBNK-BB18}, where the authors discuss the phenomenon of `phenotypic switching' of cancer cells, which can happen spontaneously, but may also be triggered by immunotherapy. 

In each of the following two sections, we first describe the new coalescent models (with spontaneous resp.\ simultaneous switching) and some of their properties, before we discuss approximating population models and their scaling limits, the seed bank diffusion (with and without jumps).

\index{simultaneous switching}
\index{responsive switching}
\index{spontaneous switching}

For completeness, we mention here that there are even more drastic ways to incorporate seed banks in population genetics. For example, a toy model extending \cite{JBNK-KKL01} to a situation with unbounded jumps was discussed in \cite{JBNK-BEGK15}. If one chooses the distribution of the jumps to be heavy-tailed with infinite variance, then a genealogy with a renewal structure emerges in the finite-expectation case. In case of an infinite expected jump size, lineages might not even coalesce at all in finite time, and the notion of a genealogy will become void. 
The mechanism of this model was discussed in \cite{JBNK-G+14} as a possible explanation for certain anomalous observations in the genome of bacterial species (`ORFan genes'), although this is debated and other explanations can be provided, see \cite{JBNK-SL18}.

\section[Seed banks with spontaneous switching]{Seed banks with spontaneous switching}\label{JBNK-sect:def}

\index{seed bank coalescent}
\index{seed bank diffusion}
\index{marked partition}

In this section, we first define the {seed bank coalescent} (with spontaneous switching) and then discuss a Wright-Fisher type population model with `strong' seed bank, corresponding to the basic migration scheme of Figure \eqref{JBNK-fig:schema}, whose ancestral process (under the usual population genetic scaling) is described by this coalescent process.
Then, we review some properties of the seed bank coalescent, in particular the time to the most recent common ancestor of a sample, and discuss the scaling limit of the frequency process of the Wright-Fisher model in the two-alleles case, the {seed bank diffusion}, which is  a classical moment dual of the seed bank coalescent. We conclude with a discussion of the long-term and boundary behaviour of this process.

\subsection{The seed bank coalescent}

The seed bank coalescent is a continuous-time Markov chain taking values in the space of partitions of natural numbers (similar to the classical Kingman coalescent), where additionally each partition block will be endowed with a `mark' from $\{a,d\}$ indicating whether the block is currently \underline{a}ctive or \underline{d}ormant. Only active blocks will be able to coalesce. Formally, for $k \ge 1$, let $\mathcal{P}_k$ be the set of partitions of $[k]:=\{1,...,k\}$. For $\pi \in \mathcal{P}_k$ let $|\pi|$ be the number of blocks of the partition $\pi.$ Define the space of \emph{marked} partitions to be 
$$
\mathcal{P}^{\{a,d\}}_k=\Big\{ (\pi, \vec{u}) \mid \pi\in \mathcal{P}_k, \vec{u} \in \{a,d\}^{|\pi|}\Big\},
$$ 
attaching to each of the $|\pi|$ blocks of a partition $\pi$ a mark from $\{a,d\}$. 

For two marked partitions $\pi,  \pi^\prime \in \mathcal{P}_k^{\{a,d\}}$, we write $\pi\succ \pi^\prime$ if $\pi^\prime$ can be constructed by merging exactly two blocks of $\pi$ carrying the $a$-mark, and the resulting block in ${\bf \pi}^\prime$ obtained from the merging both again carries an $a$-mark. For example,
$$
\big\{\{1,3\}^a\{2\}^d\{4,5\}^a\big\}\succ \big\{\{1,3,4,5\}^a\{2\}^d\big\}.
$$
We use the notation ${\bf \pi}\Join {\bf \pi}^\prime$ if ${\bf \pi}^\prime$ can be constructed by changing the mark of precisely one block of $\pi$,
 for example
$$
\big\{\{1,3\}^a\{2\}^d\{4,5\}^a\}\Join \{\{1,3\}^d\{2\}^d\{4,5\}^a\big\}. 
$$
With this notation we can now formally define the seed bank coalescent.

\begin{definition}[The seed bank coalescent]
\label{JBNK-defn:k_seedbank_coalescent}
\index{seed bank coalescent}
For $k \ge 2$ and $c,K \in (0,\infty)$ we define the \emph{seed bank $k$-coalescent} $(\Pi^{(k)}_t)_{t \ge 0}$ with \emph{switching rate} $c$ and \emph{relative seed bank size} $1/K$ to be the continuous time Markov chain with values in $\mathcal{P}_k^{\{a,d\}}$, starting in $\Pi^{(k)}_0=\{\{1\},...,\{k\}\}$ characterised by the following transitions:  
\begin{align}
\label{JBNK-eq:coalescent_transitions}
&{\bf \pi} \mapsto {\bf \pi}^\prime \,\, \text{ at rate } \,\,\begin{cases}
                                         1  & \text{ if  } {\bf \pi}\succ {\bf \pi}^\prime,\\
                                         c  & \text{ if } {\bf \pi}\Join {\bf \pi}^\prime\text{ and one $a$ is replaced by one $d$},\\
                                         cK & \text{ if } {\bf \pi}\Join {\bf \pi}^\prime\text{ and one $d$ is replaced by one $a$}.
                                        \end{cases}
\end{align}
The  \emph{seed bank coalescent} $(\Pi_t)_{t \ge 0}=(\Pi^{(\infty)}_t)_{t \ge 0}$ is then given as the
unique Markov process distributed as the projective limit as $k$ goes to infinity of the laws of the seed bank $k$-coalescents. 
\end{definition}

\begin{figure}[t]
    \centering
    \includegraphics[angle=90, width=0.85\textwidth, height=0.33\textwidth]{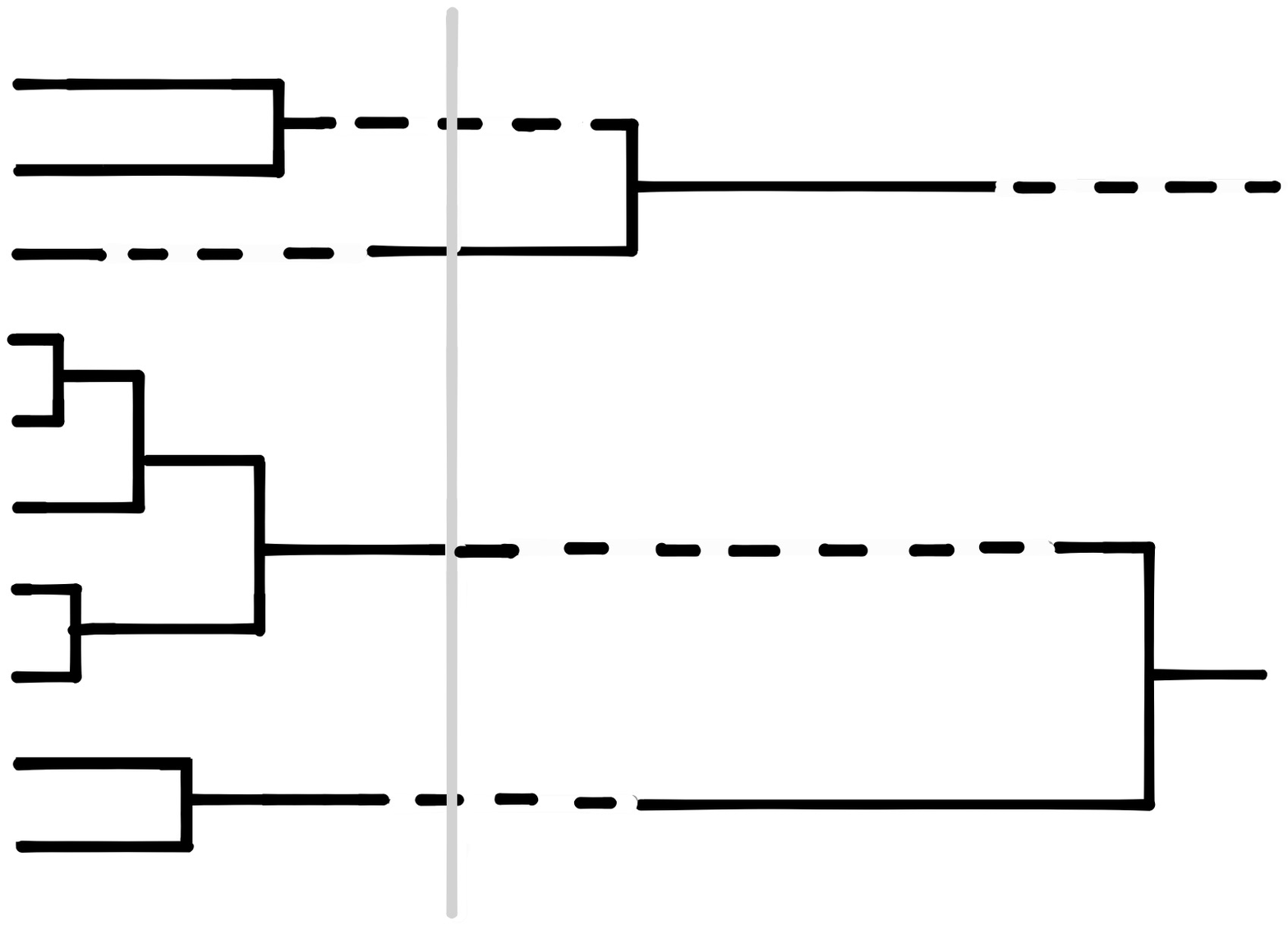}  
    \vspace{-.4cm}
    
    {\small$1 \,\,\,\,  2\,\,\,\, 3 \,\,\,\,4 \,\,\,\,5 \,\,\,\,6\,\,\,\, 7\,\,\,\, 8\,\, \,\,9\,\, \,\,10 \,\,\,$}
    \medskip
      \caption{A  realisation of the standard 10-seed bank coalescent. Dotted lines indicate `inactive lineages' (carrying a $d$-mark, which are prohibited from merging).  At the time marked with the dotted horizontal line the process is in  state $\{\{1,2\}^d\{3\}^a\{4,5,6,7,8\}^a\{9,10\}^d\}$.}
\end{figure}

The above rates state that any two \emph{active} blocks of the seed bank coalescent merge at rate 1. An active block becomes dormant at rate $c$, and a dormant block becomes active at rate $cK.$ Dormant blocks do not coalesce. This precisely describes individual switching of active lines between active and dormant states. 
The corresponding \emph{block counting process} $(N_t,M_t)_{t \ge 0}$ is the continuous time Markov chain taking values in $\NN_0\times\NN_0$ with 
transitions
\begin{equation}
\label{JBNK-eq:dual_rates}
(n,m)\mapsto \begin{cases}
(n-1,m)  & \text{ at rate } \binom{n}{2},\\
(n-1,m+1)  &\text{ at rate }  cn, \\
(n+1,m-1) & \text{ at rate }   cKm.
\end{cases}
\end{equation}

We now recall a simple Wright-Fisher type population model from \cite{JBNK-BGKW16} which is related to  Wright's two island model as mentioned in the introduction, and whose ancestral process converges under the usual scaling to the seed bank coalescent. Note that the model fits well to the basic scheme of Figure \eqref{JBNK-fig:schema} that has been described in \cite{JBNK-LJ11}.

\begin{definition}[Wright-Fisher model with strong (geometric) seed bank]\index{seed bank model}
\label{JBNK-def:seed_bank_model}
 Consider a population of constant total size $N+M$ of haploid  individuals reproducing in discrete non-overlapping generations $r=0,1,...$, consisting of a sub-population of $N$ active individuals, and a \emph{seed bank} of $M$ dormant individuals. Moreover, assume that the active population and the seed bank are of comparable size, that is, $M=M(N)=\frac{N}{K}$ for some constant $K>0$ (inserting Gauss brackets if necessary). Fix a non-negative integer $c\leq \max(N,M).$ Independently for each new generation $r+1$, individuals are obtained from those of generation $r$ by the following mechanism: 
\begin{itemize}
\item For the new active sub-population, $N-c$ active individuals are obtained by classical Wright-Fisherian (symmetric multinomial) sampling from the previous active generation. The remaining $c$ active slots are filled by sampling (without replacement) $c$ types independently and uniformly from the seed bank of the previous generation. 
\item For the new dormant sub-population, $M-c$ dormant individuals chosen uniformly at random simply persist in the seed bank, and the remaining slots are filled up by $c$ new offspring individuals sampled uniformly (with replacement) from active individuals of the previous generation.
\end{itemize}
\end{definition}

For more details on the formal description of the model see Section 1 of \cite{JBNK-BGKW16}. By construction, time spent by a single dormant individual in the seed bank is geometrically distibuted with parameter $\frac{cK}{N}.$ It is rather straightforward to check that as $N\to\infty$ (and thus $M\to\infty$ by assumption) in the usual time rescaling by the population size $N,$ the genealogy of this population converges in distribution to the seed bank coalescent (again, see \cite{JBNK-BGKW16}, Cor. 3.5 for details).  

\begin{remark}
Note that in the above model, we assumed $c$ to be an integer. To obtain the general case, that is a seed bank model and scaling limit for arbitrary $c>0$, one needs to randomize the above reproduction mechanisms. The easiest way is too choose the number of slots to be exchanged to be an independent (for each generation) binomially distributed random variable with parameters $N$ and $c/N$. This way, the expected number of slots is still $c$, and $c$ can be chosen arbitrarily from $(0, \infty)$ as $N$ grows large. We leave the details to the interested reader.
\end{remark}

Interestingly, precisely the same coalescent process as the one from Definition \ref{JBNK-defn:k_seedbank_coalescent} also arises in a peripatric speciation model from ecology investigated by Lambert and Ma in \cite{JBNK-LM15}, which the authors baptized `peripatric coalescent' in this context.
There, a dynamic metapopulation model is considered, consisting of one large founder population and several much smaller colonies (called peripheral isolates). The large population constantly produces new colonies with a certain rate, which then eventually merge back into the founder population after a certain time. Tracing ancestral lines from a suitable scaling limit of this metapopulation (at stationarity) yields the peripatric coalescent, where lines are `deactivated' when the corresponding ancestor is outside the large founder population (in some peripheral isolate), and `reactivated', when returning to the main population.

\subsection[Properties of the seed bank coalescent]{Properties of the seed bank coalescent}

Due to the fact that dormant lines do not participate in coalescence events, key quantities related to the seed bank coalescent display a qualitatively different behaviour compared to those based on the Kingman coalescent. For example, it is 
natural to expect that the \emph{time to the most recent common ancestor}
\[
T_{\text{MRCA}}:=\inf\{t>0: N_t+M_t=1\}
\]
of a sample of size $(N_0,M_0)=(n,m)$ is longer than the corresponding time for the Kingman coalescent of size $n+m$. Recall that for Kingman's coalescent, the expectation of this time bounded by 2, irrespective of the sample size. In the seed bank coalescent, if we consider a sample of $n$ active individuals, one can see that some lines may `escape' to the seed bank before participating in a merger event, where they stay inactive for some random time until they become active again and eventually (possibly after further excursions to the seed bank) merge. The $m$ dormant lines need to become active again as well, before they can merge at all. 

In fact, if the process is started with $n$ active individuals, the expected number of lines migrating to the seed bank before merging is of order $\log n.$ Time spent in the seed bank is independent for each line and exponentially distributed. Hence the amount of time it takes for the $m$ seed bank lines to activate is the maximum of $m$ independent exponentials, and thus in expectation of order $\log m.$ Combining these two observations, it is possible to give, with some additional technical effort, asymptotic bounds on the expected time to the most recent common ancestor of the seed bank coalescent, which however tend to infinity as the sample size increases.
This holds even if one starts from a sample in which all lines are active. 

\begin{theorem}[\cite{JBNK-BGKW16} Thm. 4.6]
\label{JBNK-thm:TMRCA}
Let $\EE_{(n,0)}$ be the expectation of the block-counting process $(N_t,M_t)_{t \ge 0}$ when started in $(N_0,M_0)=(n,0)$ for some $n \in \NN$. Then, 
$$
0<\liminf_{n\to\infty}\frac{\EE_{(n,0)}\big[T_{\text{MRCA}}\big]}{\log\log n}\leq \limsup_{n\to\infty}\frac{\EE_{(n,0)}\big[T_{\text{MRCA}}\big]}{\log\log n}<\infty. 
$$
\end{theorem}

\index{time to the most recent common ancestor}

In view of the above result, it is also natural to expect that the seed bank coalescent, when started from infinitely many lines, 
stays infinite for all times, and indeed we have the following result.

\begin{theorem}[\cite{JBNK-BGKW16} Thm. 4.1]
\label{JBNK-thm:cdi}
If $n+m=\infty,$ then 
$$
\PP^{(n,m)}(M_t=\infty \mbox{ for all } t \geq 0)=1.
$$
\end{theorem}

\index{coming down from infinity}

The notion of \emph{coming down from infinity} for exchangeable coalescents was introduced by Pitman \cite{JBNK-P99} and Schweinsberg \cite{JBNK-S00}, who distinguish between \emph{coming down from infinity (instantaneously)} and \emph{staying infinite (at all times)}. While the seed bank coalescent stays infinite, we will see in Section \ref{JBNK-sect:sim} that the \emph{seed bank coalescent with simultaneous switching} may stay infinite, come down from infinity instantaneously, or even come down from infinity after a finite time, depending on the choice of parameters and the initial conditions.

To asses genetic variability under the seed bank coalescent, one needs to incorporate mutations. This can be done in the standard way, by placing them 
on the active lines according to a Poisson process of rate $u/2 >0$, say, and on the dormant lines with rate $u'/2>0$. It is of course a modelling question whether mutations should be allowed in the dormant population at all, and if so, whether the rate should be reduced in comparison to the active population. To model genomic (SNP) data, one typically assumes the infinitely-many sites model.

\index{mutation rate}
\index{infinitely-many sites}

Recursions for the expected values of classical population genetic quantities such as the number of segregating sites or the number of singletons can now easily be obtained from a first step analysis, see \cite{JBNK-BEGKW15, JBNK-KR16}. From these and related quantities, various distance statistics may be calculated. Recently, a closed form representation for the site frequency spectrum for the seed bank coalescent has been obtained by Hobolth et al.\ with the help of the theory of phase-type distributions, see \cite{JBNK-HSJB18}. It turns out that for example \emph{Fu and Li's D} as well as \emph{Fay and Wu's H}, both using the full frequency spectrum, might be suitable statistics to detect the presence of seed banks \cite{JBNK-BEGKW15}, at least in the case of no or very little mutation in the dormant part of the population. 

In fact, a model selection and inference machinery, similar to the one for distinguishing among $\Lambda$- and $\Xi$-coalescents as in \cite{JBNK-BB08, JBNK-BBE13} and the contribution by Birkner and Blath in the present volume, can be derived, see \cite{JBNK-BBKWB18+} for first steps in this direction. 

\subsection{The seed bank diffusion}

For general type-space and mutation models, the scaling limit of our population model from 
Definition \ref{JBNK-def:seed_bank_model} will be a measure valued {\em seed bank Fleming-Viot process}.
However, in a bi-allelic set-up, where each individual carries one of the two possible types $\{0,1\}$, the Fleming-Viot process reduces to a two-dimensional stochastic differential equation. To describe this limit, denote by $X_r^{N}$ the frequency of the $0$-types in the active population and by $Y_r^{N}$ the frequency of the $0$-types in the dormant population in generation $r\ge 0.$ Again, it is natural to consider the limit of these frequency processes on the appropriate time-scale, and one obtains the following result:

\begin{theorem}[Seed bank diffusion, \cite{JBNK-BGKW16} Cor. 2.5] \index{seed bank diffusion} 
If $(X^N_0, Y^N_0) \to (x,y) \in [0,1]^2$ as $N \to \infty$, then the rescaled process $\big(X^N_{\lfloor Nt\rfloor}, Y^N_{\lfloor Nt\rfloor}\big)_{t\geq 0}$ converges on the Skorohod space of c\`adl\`ag paths to the unique solution of
\begin{align}
 \label{JBNK-eq:system}
 {\rm{d}} X_t & = c(Y_t -X_t) {\rm{d}}t + \sqrt{X_t(1-X_t)} {\rm{d}}B_t, \notag \\
  {\rm{d}} Y_t & = cK(X_t -Y_t) {\rm{d}}t,
 \end{align}
with $(X_0, Y_0) =(x,y) \in [0,1]^2$, where $(B_t)_{t\geq 0}$ is standard Brownian motion. 
\end{theorem}

We call $(X_t, Y_t)_{t\geq 0}$ the \emph{seed bank diffusion} with parameters $c,K>0.$ Only the active individuals reproduce and thus are subject to a Wright-Fisher noise. Activation and deactivation of individuals is governed by the rates $c$ and $cK,$ respectively, and takes the form of a classical migration term. Indeed, \eqref{JBNK-eq:system} is similar to systems of SDEs describing structured populations, in particular of \emph{island models}, \cite{JBNK-NG93, JBNK-N90, JBNK-KZH08}, a relationship which is exploited and discussed in \cite{JBNK-BBGW18}. 

A reformulation of the above system that reveals the underlying `age structure' and even provides a direct link to the model by Kaj, Krone and Lascoux \cite{JBNK-KKL01} can be stated as follows. Assume starting frequencies $X_0=x\in [0,1], Y_0=y \in [0,1]$. Then, the solution to \eqref{JBNK-eq:system} is a.s.\ equal to the unique strong solution of the stochastic delay differential equations
\begin{align*}
{\rm{d}}X_t &= c \Big( y e^{-cKt} + \int_0^t cKe^{-cK(t-s)}  X_s   {\rm{d}}s  - X_t \Big)  {\rm{d}} t + \sqrt{X_t(1-X_t)}{\rm{d}}B_t, \notag \\
{\rm{d}}Y_t &= cK \Big( -y e^{-cKt} - \int_{0}^{t} cKe^{-cK(t-s)}  X_s {\rm{d}}s + X_t \Big) {\rm{d}} t,\notag 
\end{align*}
with the same initial condition. The result rests on the fact that the driving noise is one-dimensional and can be proved via an integration-by-parts argument, see \cite{JBNK-BBGW18}, Prop. 1.4. The second component is now just a deterministic function of the first, which in integral form reads
$$
Y_t= ye^{-cKt} + \int_0^t cK e^{-cK(t-s)}\, X_s {\rm{d}}s,
$$ 
and the first delay equation is independent of the second. The delay representation allows an elegant interpretation of the ``time-lag'' caused by a seed bank.
Indeed, in integral form it is given by 
$$
X_t = x+\frac yK (1-e^{-cKt})- c\int_0^t  e^{-cK(t-s)}\, X_s {\rm{d}}s + \int_0^t \sqrt{ X_s(1-X_s)} \, {\rm{d}}B_s,
$$
and shows that the genetic type of any currently reactivated individual is determined by the corresponding type frequency of the active population alive at an exponentially distributed time ago (cut off at time 0). This is in line with the fact that the time that individuals should spend in the seed bank is given by an exponentially distributed random variable with parameter $cK$, which is the scaling limit of the geometric seed bank age distribution (with parameter $cK/N$) in the population model on the new evolutionary time-scale measuring time in units of order $N$. 

\medskip

The seed bank diffusion $(X_t, Y_t)_{t\geq 0}$ is the {\em moment dual} of the block counting process $(N_t, M_t)_{t\geq 0}$ in a classical sense.

\index{seed bank diffusion}
\index{moment duality}
\index{duality}

\begin{proposition}[\cite{JBNK-BGKW16} Thm.\ 2.8]
\label{JBNK-eq:dual}
For every $(x,y)\in [0,1]^2 $, every $n,m\in \NN_0$ and every $t\geq 0$
\begin{equation}\label{JBNK-eq:duality}
\mathbb{E}_{x,y}\big[X_t^n Y_t^m\big]=\mathbb{E}^{n,m}\big[x^{N_t} y^{M_t}\big],
\end{equation}
where $\mathbb{E}_{x,y}$ denotes the expectation with respect to the measure $\mathbb{P}_{x,y}$ for $(X_t,Y_t)_{t\geq 0}$ started at $X_0=x, Y_0=y,$ and  $\mathbb{E}^{n,m}$ refers to $(N_t, M_t)_{t\geq 0}$ started in $(n,m).$
\end{proposition}

Moment duality is a useful tool to study the long-term behaviour of population models. 
A classical Wright-Fisher diffusion starting in some value $z \in [0,1]$ will get absorbed at the boundaries after finite time a.s.\  in fact with finite expectation, hitting 1 with probability $z$. The situation is slightly more involved for our frequency process in the presence of a strong seed bank. Obviously, $(0,0)$ and $(1,1)$ are the only absorbing states for the system.  

Let us first consider the long-term behaviour in law.  Applying moment duality \eqref{JBNK-eq:duality}, one finds that all mixed moments of $(X_t,Y_t)_{t \ge 0}$ converge to the \emph{same} finite limit depending only on $x,y, K$. More precisely, for each fixed $n,m\in\NN$, 
\begin{equation}
\label{JBNK-eq:moment_value}
\lim_{t \to \infty} \mathbb{E}_{x,y}[X_t^{n}Y_t^{m}] = \lim_{t \to \infty} \mathbb{E}^{n,m}\big[x^{N_t} y^{M_t}\big] = \frac{y+xK}{1+K},
\end{equation}
since the block counting process, when started with finitely many individuals, will always collapse to the situtation in which only one line remains that switches among the active and dormant state cf.\ \cite{JBNK-BGKW16}, Section 2.3 for details. From this, fixation in law follows by uniqueness of the moment problem on $[0,1]^2$ and the Stone-Weierstra\ss\ Theorem and we get

\begin{corollary}[\cite{JBNK-BGKW16} Cor. 2.10]
\label{JBNK-cor:fix_law}
Given $c, K$, $(X_t, Y_t)$ converges in distribution as $t\to\infty$ to a two-dimensional random variable $(X_\infty, Y_\infty),$ whose distribution is given by
\begin{equation}
\label{JBNK-eq:momentconvergence}
 \mathcal{L}_{(x,y)}\big( X_\infty, Y_\infty \big) =  \frac{y+xK}{1+K} \delta_{(1,1)} + 
\frac{1+(1-x)K-y}{1+K} \delta_{(0,0)}.
\end{equation}
\end{corollary}
This is in line with the classical results for the Wright-Fisher diffusion: As $K \to \infty$ (that is, the seed bank becomes small compared to the active population), the fixation probability of
$a$ alleles approaches $x$. Further, for small $K$ (that is, large seed bank), the fixation probability is governed by the initial fraction $y$ of $a$-alleles in the seed bank.

Observing that 
$$
KX_t+Y_t =Kx+y+K\int_0^t \sqrt{X_s(1-X_s)} \, {\rm{d}}B_s, \quad t \ge 0,
$$ 
gives rise to a bounded martingale, and given the shape of the limiting law \eqref{JBNK-eq:momentconvergence}, one can also get almost sure convergence of $(X_t, Y_t)$ to $(X_\infty, Y_\infty)$ as $t\to\infty.$ However, as we will see as a special case of Theorem \ref{JBNK-thm:boundary} below, fixation will not happen in finite time. The intuition behind this can also be directly understood from \eqref{JBNK-eq:system}, where we can compare the seed-component $(Y_t)_{t\geq 0}$ to the solution of the deterministic equation
\[
{\rm{d}}y_t=-cKy_t{\rm{d}}t,
\]
corresponding to a situation where the drift towards 0 is maximal (or to ${\rm{d}}y_t=cK(1-y_t){\rm{d}}t$ where the drift towards 1 is maximal). Since $(y_t)_{t\geq 0}$ does not reach 0 in finite time if $y_0>0,$ neither does $(Y_t)_{t\geq 0}.$ This is also reflected in the fact that the
block-counting process $(N_t,M_t)_{t \ge 0}$, started from an infinite initial state, \emph{does not come down from infinity}.

\medskip

We now introduce mutation and keep our focus again on the two alleles model. We assume that in the active population, mutation from type $0$ to type $1$ happens at rate $u_1,$ and from $1$ to $0$ at rate $u_2.$ The respective rates in the dormant population are denoted by $u_1', u_2',$ which may be different from the active population (or even vanish). 
Then, we obtain the following system of SDEs.

\begin{definition}[Seed bank diffusion with mutation]
\label{JBNK-defn:system-mut}
 The {\em seed bank diffusion with mutation} with parameters $u_1,u_2,u_1',u_2',c,K$ is given by the unique strong solution of the initial value problem 
\begin{align}
\label{JBNK-eq:system-mut}
{\rm{d}} X_t & = \big[-u_1X_t +u_2(1-X_t)+ c(Y_t -X_t)\big]{\rm{d}}t + \sqrt{X_t(1-X_t)}{\rm{d}}B_t, \notag \\[.1cm]
{\rm{d}} Y_t & = \big[-u_1'Y_t+u_2'(1-Y_t) + cK(X_t -Y_t)\big]{\rm{d}}t,
\end{align}
with $(X_0, Y_0) =(x,y) \in [0,1]^2$. 
\end{definition}

Note that in order to study the boundary behaviour of the seed bank diffusion with mutation, one cannot simply refer to Feller's boundary classification machinery, resting on speed measure and scale function, since the above system is two-dimensional. Still, it is possible to provide a rather satisfactory characterization. To this end, define the first hitting time of $X$ of the boundary 0 by
\begin{align*}
 \tau^X_0:=\inf\{t \geq 0 \mid X_t=0\},
\end{align*}
and define $\tau^X_1$, $\tau^Y_0$ and $\tau^Y_1$ analogously. We say that \emph{$X$ will never hit 0} (from the interior), if for every initial distribution $\mu_0$ such that $\mu_0((0,1)^2)=1$, we have
\begin{align*}
  \PP^{\mu_0}\left( \tau^X_0 < \infty\right) = 0.
\end{align*}
Using similar notation for the other cases, the following boundary classification can be achieved.

\begin{theorem}\label{JBNK-thm:boundary}
 Let $(X_t,Y_t)_{t \geq 0}$ be the solution to \eqref{JBNK-eq:system-mut} with parameters satisfying $u_1,u_2,u_1',u_2'\geq 0$ and $c,K>0$. 
 \begin{enumerate}
  \item Started from the interior $X$ will never hit 0 if and only if $2u_2 \geq 1$.
  \item Started from the interior $X$ will never hit 1 if and only if $2u_1 \geq 1$.
  \item Started from the interior $Y$ will never hit 0.
  \item Started from the interior $Y$ will never hit 1.
 \end{enumerate}
\end{theorem}
For a more general result, including also the two-island diffusion case, and clarifying the roles of the different parameters in more detail, see Theorem 3.1 in \cite{JBNK-BBGW18}, which uses recent progress in the theory of polynomial diffusions \cite{JBNK-FL16} and a version of ``McKean's argument''. The latter is particularly useful since it is based on submartingale convergence arguments that work in the multi-dimensional case (as opposed to the speed-measure scale-function formalism which is restricted to dimension one) and seems to have appeared first in \cite[Problem 7, p.47]{JBNK-M69}.

\section[Simultaneous switching]{Seed banks with simultaneous switching}
\label{JBNK-sect:sim}

\subsection[Simultaneous switching]{The seed bank coalescent  with simultaneous switching}

We now extend the seed bank coalescent to incorporate simultaneous switching (cf.\ \cite{JBNK-BGKW19}), as discussed in the introduction. In addition to the parameters $c$ and $K$ from Section \ref{JBNK-sect:def}, we now also fix two finite measures $\Lambda$ and $\overline{\Lambda}$ on $[0,1]$ which govern the sizes of \emph{simultaneous switching events} from active to dormant, and vice versa. To keep the analogy with migration models, we also call these {\em large migration events}.

Again our seed bank coalescent with simultaneous switching will be a continuous time Markov chain on the space of 
marked partitions $\mathcal{P}^{\{a,d\}}$. In addition to the notation in the previous Section \ref{JBNK-sect:def}, we write ${\bf \pi}\Join_m {\bf \pi}^\prime$ if ${\bf \pi}^\prime$ can be constructed by changing the mark of precisely $m$ blocks of $\pi$ from $a$ to $d,$ and ${\bf \pi}\Join^l {\bf \pi}^\prime$ if ${\bf \pi}^\prime$ can be constructed by changing the mark of precisely $l$ blocks of $\pi$ from $d$ to $a$. 

\begin{definition}[The seed bank coalescent with simultaneous switching]\index{simultaneous switching}
\label{JBNK-defn:seedbank_coalescent-sim}
Fix $c, K \in (0,\infty)$ and finite measures $\Lambda,\overline{\Lambda}$ on $[0,1]$ such that $\Lambda(\{0\})=\overline{\Lambda}(\{0\})=0.$ For $k\geq 1$ we define the \emph{seed bank $k$-coalescent with simultaneous switching} to be the continuous time Markov chain with values in the space of marked partitions $\mathcal{P}_k^{\{a,d\}}$, characterised by the following transitions:  
\begin{align}
\label{JBNK-eq:coalescent_transitions-sim}
&{\bf \pi} \mapsto {\bf \pi}^\prime \text{ at rate } \begin{cases}
                                         1  & \text{ if  } \pi \succ \pi',
                                         \\
                                         c+\int_{[0,1]}z(1-z)^{{|\pi|-1}}\frac{\Lambda(\dd z) }{z} & \text{ if } \pi\Join_1 \pi',\\
                                         c K+\int_{[0,1]}z(1-z)^{{|\pi|-1}}\frac{\overline{\Lambda}(\dd z)}{z} & \text{ if } \pi\Join^1 \pi',\\
                                         \int_{[0,1]}z^k(1-z)^{|\pi|-k} \frac{\Lambda(\dd z)}{z} & \text{ if } \pi\Join_k \pi', 2\leq k\leq |\pi|, \\
                                         \int_{[0,1]}z^l(1-z)^{|\pi|-l}\frac{\overline{\Lambda}(\dd z)}{z}  & \text{ if } \pi\Join^l \pi', 2\leq l\leq |\pi|.\\
                                        \end{cases}
\end{align}
The  \emph{seed bank coalescent with simultaneous switching} is given by the projective limit of seed bank $k$-coalescents with simultaneous switching as $k \to \infty$.
\end{definition}

Note that we recover the ordinary seed bank coalescent for the choice $\Lambda=\overline{\Lambda}\equiv 0.$ The corresponding block counting process then  has transitions
\begin{equation}
\label{eq:dual_rates}
(n,m)\mapsto \begin{cases}
(n-1,m)  & \text{ at rate } \binom{n}{2}, \;  n\geq 2\\
(n-1,m+1)  &\text{ at rate }  \left(c+\int_0^1z(1-z)^{n-1}\frac{{\Lambda}(\dd z)}{z} \right)n, \; n\geq 1 \\
(n-k,m+k)  &\text{ at rate }  \binom{n}{k}\int_0^1 z^k(1-z)^{n-k}\frac{{\Lambda}(\dd z)}{z} ,\;  2\leq k\leq n, \\
(n+1,m-1) & \text{ at rate }  \left (c K +\int_0^1z(1-z)^{m-1}\frac{\overline{\Lambda}(\dd z)}{z} \right)m, \;  m\geq 1\\
(n+l,m-l)  &\text{ at rate }  \binom{m}{l}\int_0^1z^l(1-z)^{m-k}\frac{\overline{\Lambda}(\dd z)}{z}, \; 2\leq l\leq m. 
\end{cases}
\end{equation}

\begin{figure}[t]
\begin{center}
\includegraphics[scale=0.3]{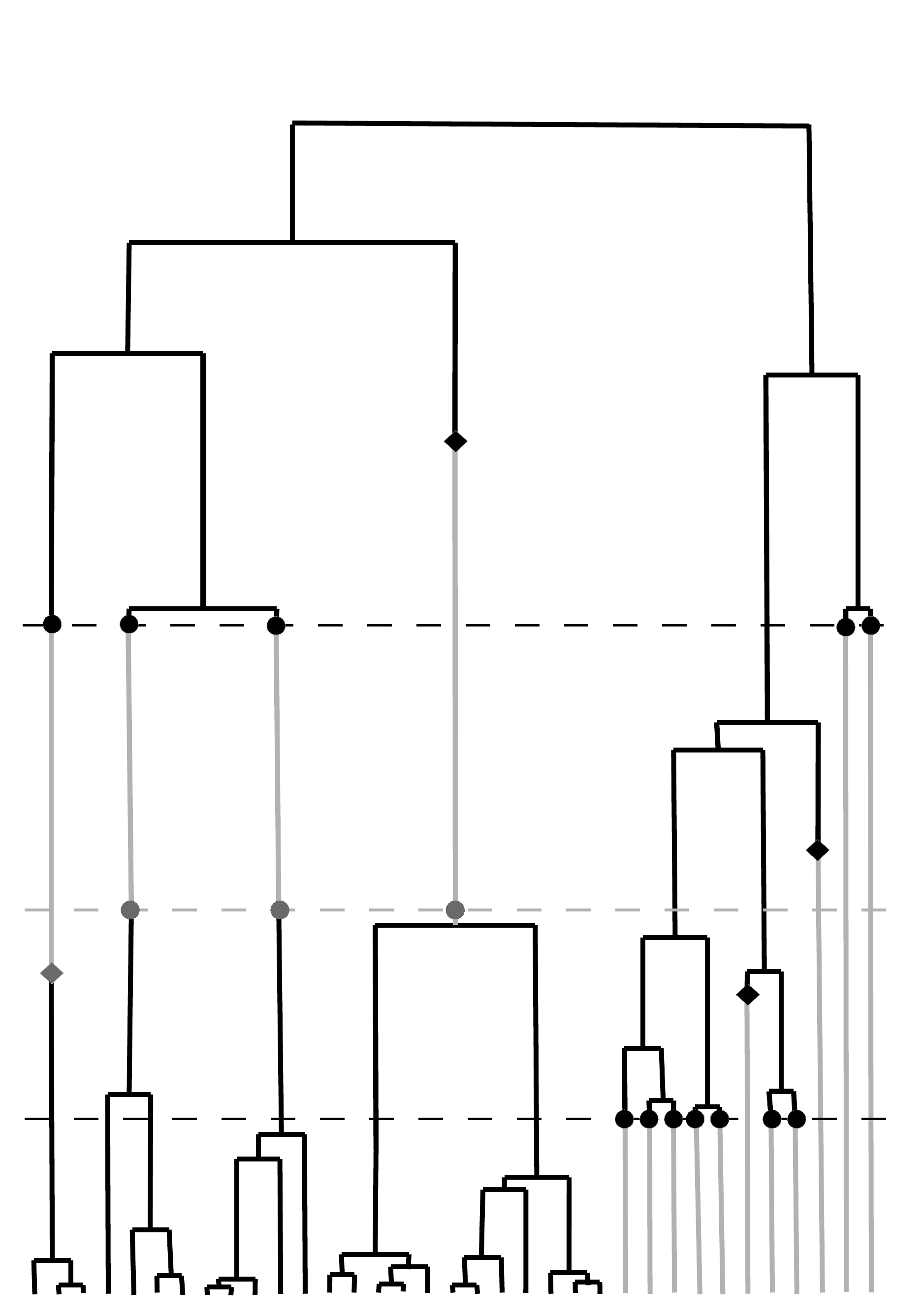}
\caption{\label{JBNK-fig:coalsim} Seed bank coalescent with simultaneous switching. Dormant lines in grey. Dots indicate switches, dashed lines simultaneous switching events.
} 
\end{center}
\end{figure}

The process is thus an extension of the seed bank coalescent defined in Definition \ref{JBNK-defn:k_seedbank_coalescent}, where in addition to the spontaneous switching of single lines at rate $c$ resp.\ $cK$ there are coordinated switches of a (large) number of lines. In an event of the latter type, a number $z\in [0,1]$ is determined according to the measure $z^{-1}\Lambda(\dd z)$ resp. $z^{-1}\bar{\Lambda}(\dd z)$, and then each of the $n$ resp. $m$ lines determines independently with probability $z$ whether or not to participate in the switching, leading to a binomial number of lines changing state. Figure \ref{JBNK-fig:coalsim} shows an instance of such a coalescent.

\begin{remark}
Note that the (possibly infinite) measure $z^{-1}\Lambda(\dd z)$ is formally similar to the measure $z^{-2}\Lambda(\dd z)$ driving the jumps
of a $\Lambda$-coalescent (see, e.g.\ \cite{JBNK-P99, JBNK-S99, JBNK-DK99}; cf.\ also the article by Birkner and Blath in the present volume), but the singularity at 0 is at most of order 1 instead of order 2 as for multiple merger events. The intuitive reason for this is that coalescence events always require to involve at least two particles in order to be visible, whereas migration events are visible already if only one particle is affected. 
\end{remark}

The seed bank coalescent with simultaneous switching is again the limiting genealogy of a Wright-Fisher type population model (see \cite{JBNK-BGKW19} for details). In addition to $N,M,c,K$ from the previous section, fix probability measures $\mu_N, \overline{\mu}_N$ on $[0,1]$ (which will later be linked to $\Lambda$ and $\bLambda$ of Definition \ref{JBNK-defn:seedbank_coalescent-sim}). In each generation one of the following events takes place (independent between generations):
\begin{itemize}
\item[\textbf{S}] \emph{Small-scale migration event} (of size $o(N)$) between active and dormant according to the mechanism for the ``ordinary'' seed bank coalescent.
\item[\textbf{F}] \emph{Simultaneous switching from dormant to active:}
Sample $z\in [0,1]$ according to a probability measure $\mu_N(\dd z).$ For the new active generation, $(1-z)N$ active individuals are obtained by multinomial sampling from the previous active generation.
The remaining $zN$ active slots are filled by sampling (without replacement) $zN$ types independently and uniformly from the seed bank types of the previous generation. The seed bank stays as it is.
\item[\textbf{D}] \emph{Simultaneous switching from active to dormant:}
Sample $z\in [0,1]$ according to a probability measure $\overline{\mu}_N(\dd z).$ The $N$ active individuals in the next generation are produced by multinomial sampling from the active individuals in the previous generation. 
For the new seed bank generation, $z M$ dormant individuals from the previous generation are replaced by new dormant individuals obtained by multinomial sampling from the previous active generation. The remaining $(1-z)M$ dormant individuals stay in the seed bank.
\end{itemize}

\smallskip

In order to obtain a non-trivial limit of the above seed bank model with simultaneous switching, we need to make some scaling assumptions. In particular, the large migration events are required to happen much more rarely than the small migration events. Let $(r_N)_{N\in\NN}$ and  $(\br_N)_{N\in\NN}$ denote sequences of non-negative numbers such that $(r_N/N)_{N\in\NN}$ and $(\br_N/N)_{N\in\NN}$ converge to 0 as $N\to\infty$ (recall the assumption $M=N/K$ from Definition \ref{JBNK-def:seed_bank_model} when taking the limits). We assume that in each generation, an event of type $F$ happens with probability $r_N/N,$ and an event $D$ with probability $\br_N/N.$ Hence, events of type $S$ happen with probability $1-(r_N+\br_N)/N.$ Assume moreover that we have the weak limits
\begin{equation}
\lim_{N\to\infty}r_N\mu_N(\dd z)=z^{-1}\Lambda(\dd z)
\end{equation}
and analogously for $\br_N, \bmu_N, \bLambda.$ Then, the ancestral process of the above population model will converge weakly to a seed bank coalescent with simultaneous switching and switching measures $\Lambda$ and $\bLambda$. Note that the (total mass of) $\bLambda$ depends on $K.$

\smallskip

An analogous result holds forward in time, giving rise to a {\em seed bank Fleming-Viot process with simultaneous switching} (in the general type space case), which boils down to a seed bank diffusion with jumps in the two alleles case that we describe now.

\begin{theorem}[Seed bank diffusion with simultaneous switching, \cite{JBNK-BGKW19} Thm.\ 1.6]
\label{JBNK-thm:conv_random}
Under the above assumptions, the rescaled frequency process $\big(X^N_{\lfloor Nt\rfloor}, Y^N_{\lfloor Nt\rfloor}\big)_{t\geq 0}$ converges on the Skorohod space of c\`adl\`ag paths to the unique solution of
\begin{align}
\label{JBNK-eq:system_variable}
{\rm{d}} X_t & = c(Y_t -X_t)\, {\rm{d}}t + \sqrt{X_t(1-X_t)}\, {\rm{d}}B_t\,  \\[.1cm]
	          & \qquad \qquad\qquad\qquad \quad \, + \int \limits_{ [0,1]} z\big(Y_{t-}-X_{t-}\big) \, N^{F}\big({\rm{d}}t,{\rm{d}}z \big),\notag \\[.1cm]
{\rm{d}} Y_t & = c K(X_t -Y_t)\, {\rm{d}}t + \int \limits_{[0,1]} z\big(X_{t-}-Y_{t-}\big) \, 
   N^{D}\big({\rm{d}}t,{\rm{d}}z\big), \notag 
\end{align}
with $(X_0, Y_0) =(x,y) \in [0,1]^2$, where $(B_t)_{t\geq 0}$ is a standard Brownian motion and $(N^{F}(r, t))_{t\geq 0}$ and $(N^{D}(r,t))_{t\geq 0}$ are independent standard Poisson point processes on $(0, \infty) \times [0,1]$ with intensity measures
$$
\lambda(\dd t) \otimes z^{-1}\Lambda(\dd z) \quad \mbox{ resp.\ } 
\quad
\lambda(\dd t) \otimes z^{-1} \bar{\Lambda}(\dd \bz).
$$ 
Here, $\lambda$ denotes the Lebesgue measure on $\RR.$ The integrals in \eqref{JBNK-eq:system_variable} are taken with respect to ${\rm{d}}z$. Moreover, $(X_t, Y_t)_{t\geq 0}$ is the moment dual of $(N_t, M_t)_{t\geq 0}$, that is, the processes satisfy the analogue of \eqref{JBNK-eq:duality}.
\end{theorem}

In a similar manner as for the seed bank coalescent with spontaneous switching one obtains fixation probabilities by applying moment duality. The almost sure boundary behaviour of the seed bank diffusion with jumps seems largely open.

\smallskip

We conclude with a result on the tree properties of the seed bank coalescent with simultaneous jumps, which are less well analysed than those of the seed bank coalescent, but seem to exhibit a more complex structure. 

Indeed, in Theorem \ref{JBNK-thm:cdi} we saw that seed bank coalescent does not come down from infinity (neither instantaneously nor after a finite time), due to the fact that even within a very short time, infinitely many active lines may escape to the seed bank, from where it takes long to come back. It turns out that in the case of simultaneous switching, there is a regime with qualitatively different behaviour. 

\begin{theorem}[\cite{JBNK-BGKW19} Theorem 2.7]
\label{JBNK-thm:MRcomingdown}
Consider the block-counting process $(N_t, M_t)_{t \ge 0}$ of the seed bank coalescent with simultaneous switching under the assumptions of Definition \ref{JBNK-defn:seedbank_coalescent-sim}. Let $Y$ be a random variable with distribution $\frac{1}{\Lambda([0,1])}\Lambda$. 
\begin{itemize}
\item[(a)] If $\bLambda(\{1\})=0$, then the block-counting process started in $(N_0, M_0)=(n,\infty), n\in \NN_0\cup\{\infty\}$ will stay infinite for all times.
\item[(b)] If the block-counting process is started in $(\infty, m), m\in\NN_0,$ then the process comes down from infinity instantaneously if $\EE[-\log(Y)] < \infty$ and $c=0$. If $\EE[-\log(Y)] =\infty$ or $c>0,$ it stays infinite for all times.
\item[(c)] If $\bLambda(\{1\})>0, c=0$ and $\EE[-\log(Y)] < \infty,$ then the block-counting process started from $(n,\infty), n\in\NN_0\cup \{\infty\}$ comes down from infinity after a finite time, but not instantaneously. 
\end{itemize} 
\end{theorem}

Part (a) of this theorem just states that the block-counting process started with infinitely many dormant lines stays infinite unless there is the possibility of emptying the seed bank at once. This is due to the fact that otherwise the overall migration rate from dormant to active is linear, and thus the seed bank will always stay infinite. If started with finitely many seeds only, then it depends on the specifics of $c, \Lambda$ and $\bLambda$ whether the process stays infinite or comes down from infinity (instantaneously or after a finite time). Note that $c=0$ is a necessary condition for coming down from infinity, otherwise a simple comparison with the seed bank coalescent without simultaneous switching shows that the process stays infinite. The condition $\EE[-\log Y]<\infty$ and the proof of the theorem are inspired by a similar condition of Griffiths \cite{JBNK-G14} for Lambda-coalescents. It ensures that the switching rates from active to dormant are sufficiently small compared to the total transition rates, such that a Borel Cantelli argument ensures that there are only finitely many transitions from active to dormant before the remaining active lines have coalesced.

\section{Open problems and perspectives for future work}

In this section we briefly mention some possible model extensions and open problems related to seed banks.

First, note that the additional feature that the lineages of a coalescent process can be either active or dormant can certainly be extended from the Kingman-coalescent framework to the much more general class of Xi-coalescents (cf.\ \cite{JBNK-S2000} or the contributions of Birkner \& Blath and Kersting \& Wakolbinger in this volume. Together with Lambert and Ma who defined the `peripatric coalescent' in \cite{JBNK-LM15}, we suggest to call this new class `on/off coalescents', since the possibility to take part in a coalescence event of a lineage can be turned `on' (active line) and `off' (dormant line). Such an {on/off coalescent} will then be defined by a pair of parameters $c, c'>0$ (where $c'=cK$ in the seed bank coalescent) describing the spontaneous switching rates, a finite measure $\Xi$ on the infinite-dimensional simplex, describing the simultaneous multiple collisions, and two finite measures $\Lambda, \bLambda$ on $[0,1]$ with no atom at 0, describing the simultaneous switching events. An investigation of the properties of this class of processes can provide a rich source of problems to probabilists with an interest in coalescent processes.

\index{on/off coalescent}
\index{peripatric coalescent}

Another active line of research is to investigate seed bank models in a spatial context, for example on the discrete torus or on the hierarchical group, see the contribution of Greven \& den Hollander in this volume. In the continuum, a natural idea is to incorporate a seed bank into the Fisher-KPP equation from ecology, whose dual is given by branching Brownian motions. We expect that the dual of a `Fisher-KPP equation with dormancy' should give rise to an `on-off'-branching Brownian motion, and this feature should affect the dynamic behaviour of the solutions of the equation. For example, we expect that the wave-speed of travelling wave solutions (if they still exist) should be significantly reduced. This is currently research in progress.

A different route to follow would be to incorporate seed banks/dormancy into the rapidly developing framework of adaptive dynamics models, see the contribution of Bovier in this volume. A very concrete question here would be to investigate under which conditions a newly emerged `dormancy-trait' could fixate in a population under competitive and selective pressure, in particular when the dormancy-trait comes with an evolutionary cost (such as a significantly reduced reproductive fitness). Again, this is research in progress.

Many further interactions of dormancy with other evolutionary forces are possible. In fact, as mentioned in the introduction, the effect of dormancy affects the macroscopic behaviour of populations in so many different ways that one is tempted to consider this effect an evolutionary force itself.

\medskip

{\bf Acknowledgement.} The authors are grateful for the comments and remarks of two anonymous referees and the editors, which significantly improved this manuscript.

 \end{document}